\documentclass[12pt]{amsart}

\usepackage{amsmath,amsthm,amsfonts,amssymb,mathtools,mathrsfs,bm}
\usepackage{amscd}
\usepackage{enumerate}
\usepackage{xcolor}
\definecolor{colordelink}{rgb}{0,0,0.50}
\definecolor{colordecite}{rgb}{0,0.5,0}
\definecolor{colordeurl}{rgb}{0,0.41,0.5}
\usepackage[pagebackref=true]{hyperref}
\hypersetup{
    colorlinks=true,
    linkcolor=colordelink,
    citecolor=colordecite,
		urlcolor=colordeurl,
}
\usepackage[capitalize,nameinlink,noabbrev]{cleveref}
\usepackage[all]{hypcap}

\usepackage[utf8]{inputenc}

\usepackage[T1]{fontenc}
\usepackage{lmodern}

\usepackage{todonotes}

\usepackage{tikz}
\usetikzlibrary{decorations.markings,arrows.meta,calc}
\usepackage{verbatim}

\def\dim{\operatorname{dim}}

\newcommand{\RR}{\mathbb{R}}

\newcommand{\CC}{\mathbb{C}}

\newcommand{\C}{\mathbb{C}}

\renewcommand{\Bbb}{\mathbb{B}}

\newcommand{\Dbb}{\mathbb{D}}

\newcommand{\Gcal}{\mathcal{G}}
\newcommand{\Hcal}{\mathcal{H}}

\newcommand{\Scal}{\mathcal{S}}

\newcommand{\Ucal}{\mathcal{U}}

\makeatletter
\newcommand{\tpitchfork}{%
  \vbox{
    \baselineskip\z@skip
    \lineskip-.52ex
    \lineskiplimit\maxdimen
    \m@th
    \ialign{##\crcr\hidewidth\smash{$-$}\hidewidth\crcr$\pitchfork$\crcr}
  }%
}
\makeatother

\theoremstyle{plain}
\newtheorem{theorem}{Theorem}[section]
\newtheorem*{theorem*}{Theorem}
\newtheorem{lemma}[theorem]{Lemma}
\newtheorem*{lemma*}{Lemma}

\newtheorem*{corollary*}{Corollary}
\newtheorem{proposition}[theorem]{Proposition}

\newtheorem{theo}{Theorem}[section]

\newtheorem{cor}[theorem]{Corollary}
\newtheorem{prop}[theorem]{Proposition}

\theoremstyle{definition}
\newtheorem{definition}[theorem]{Definition}
\newtheorem{proposition/definition}[theorem]{Proposition/Definition}

\newtheorem{example}[theorem]{Example}

\newtheorem{remark}[theorem]{Remark}

\Crefname{notation}{Notation}{Notations}

\theoremstyle{remark}

\def\e{{\epsilon}}

\newcommand{\norm}[1]{\lVert #1\rVert}

\begin{document}
\title[]{Embedded Topological Triviality of Separable Families of Singularities}
\author{R. Giménez Conejero}
\author{A. Lind}
\author{A. Menegon}

\address{Roberto Giménez Conejero: Department of Engineering, Mathematics and Science Education - Mid Sweden University.}
\email{roberto.gimenez@uv.es}
\address{Andreas Lind: Department of Engineering, Mathematics and Science Education - Mid Sweden University.}
\email{andreas.lind@miun.se}
\address{Aur\'elio Menegon:
Department of Engineering, Mathematics and Science Education,
Mid Sweden University, Sundsvall, Sweden;
and Departamento de Matem\'atica,
Universidade Federal da Para\'iba,
Jo\~ao Pessoa, Brazil.}
\email{aurelio.menegon@miun.se}


\begin{abstract}
Understanding how singularities behave under small perturbations is a central theme in singularity theory.
In this paper, we study separable families $f_s=f_0+\rho(s)\varphi$ of analytic function-germs on a germ of a complex analytic space $(X,0)$, and we give sufficient conditions for them to admit a topological trivialization.
Our results extend previous work by the third author and collaborators, moving from abstract triviality to the embedded setting.
As an application, we obtain new instances of topological stability, including a class of $\mu$-constant deformations.
These findings provide a new perspective on the long-standing $\mu$-constant conjecture, one of the major open problems in the field.
\end{abstract}

\maketitle

\section{Introduction}
\label{section_0}

Let $(X,0)$ be a germ of a complex analytic space in $\C^N$, with $\dim X >1$, and consider a family of analytic function-germs
\[
f_s: (X,0) \to (\C,0),
\]
depending analytically on a parameter $s \in \C$. The corresponding zero sets
\[
V_s \coloneqq f_s^{-1}(0)
\]
define a family of analytic germs $\big(V_s,0\big) \subset (X,0)$. A central question in singularity theory is to determine conditions under which the germs $\big(V_s,0\big)$ are topologically equivalent for small values of $s$, i.e., there exists a germ of homeomorphism $h:(X,0)\to(X,0)$ such that $(V_0,0)$ is mapped to $(V_s,0)$.

In the classical case where $X=\C^n$, each $f_s$ has an isolated singularity, and $n \neq 3$, L\^e and Ramanujam proved in their celebrated paper \cite{Le1976} that the family $(f_s)$ has constant topological type if and only if the Milnor number $\mu(f_s)$ remains constant. Later, Timourian proved in \cite{Timourian1977} that these are, in turn, equivalent to having topological triviality in the family. Their proof, however, relies on the $h$-cobordism theorem and therefore does not apply in dimension three. Since then, many partial results have been obtained (see for instance \cite{Perron1999,Parusinski1999,Borodzik2014,LeytonAlvarez2022}), but the general problem remains open. This is the well-known open problem called the $\mu$\textit{-constant conjecture} or the \textit{Lê-Ramanujam conjecture}, which has become one of the most challenging and fascinating open problems in the field. Very recently, a manuscript generated by an AI model \cite{OpenAI2026} claims a proof of the remaining case $n=3$, namely, that every one-parameter family of isolated hypersurface singularities in $\C^3$ with constant Milnor number is topologically trivial. We stress that this manuscript has not been peer reviewed and, to the best of our knowledge, it has not yet been checked by human experts.

For more general analytic spaces $(X,0)$, the problem is even harder, and most things are unsolved. In order to make progress, one usually restricts attention to special classes of deformations. A natural case to consider is that of linear deformations
\[
F(z,s) = f_0(z) + s \varphi(z),
\]
where $\varphi:(X,0)\to(\C,0)$ is an analytic function-germ and $f_s(z) = F(z,s)$. In this setting, it was proved in \cite{Menegon2023} that if $(f_s)$ has a uniform critical set $\Sigma=\{0\}$ with respect to some $(w)$-regular stratification $\mathcal S$ of $X$ (in the sense of \cite{Verdier1976}), then the family has constant \emph{abstract} topological type. That is, there exists a neighborhood $W\subset\C$ of $0$ and homeomorphisms of germs
\[
h_s:\,(V_s,0) \to (V_0,0), \quad s\in W.
\]

This extends Parusi\'nski’s theorem \cite{Parusinski1999} (and hence the result of L\^e–Ramanujam \cite{Le1976}) to more general analytic sets $X$. On the other hand, Parusi\'nski’s and L\^e–Ramanujam’s theorems are stronger in the sense that they provide an \emph{embedded} trivialization: an ambient homeomorphism of $(\C^N,0)$ sending $(V_s,0)$ onto $(V_0,0)$.

The aim of this paper is to extend the techniques of \cite{Menegon2023} towards the construction of embedded trivializations, thereby bridging the gap with the classical results. We state our results for \emph{separable families} of the form
\[
F(z,s) = f_0(z) + \rho(s)\,\varphi(z),
\]
where $\rho:(\C,0)\to(\C,0)$ and $\varphi:(X,0)\to(\C,0)$ are analytic function-germs. Linear families appear as the special case $\rho(s)=s$. Conversely, a separable family is obtained from the linear family $f_0+t\,\varphi$ by the change of parameter $t=\rho(s)$, so we prove our results for linear families and deduce the separable case from them (see \cref{section_5}).

Our results are formulated in terms of the map
\[
G=(f_0,\varphi):(X,0)\to(\C^2,0)
\]
and its discriminant $\Delta_G$ (see \cref{section_2}). Indeed, $V_s=G^{-1}(H_{\rho(s)})$, where $H_t=\{y_0+t\,y_1=0\}\subset\C^2$, so the family $(V_s)_s$ is obtained by pulling back the family of lines $(H_{\rho(s)})_s$ by $G$. Set
\[
N\coloneqq G^{-1}(\Delta_G).
\]
Our first result trivializes the family, together with its ambient space, outside $N$ (see also \cref{rem_baselocus}).

\begin{theo}\label{theo_A}
Let $f_s:(X,0)\to(\C,0)$ be a separable family of analytic function-germs which is $\Delta$-regular with respect to a Whitney stratification $\Scal$ of $X$ (for instance, if $(f_s)$ has a uniform critical set), and assume that $G$ is transverse to small spheres off $\Delta_G$ (see \cref{def:ML}). Then the family is topologically trivial outside $N$: there exist a neighborhood $M$ of $0$ in $X$, a neighborhood $D$ of $0$ in $\C$ and a homeomorphism
\[
(M\setminus N)\times D\to(M\setminus N)\times D,\qquad (x,s)\mapsto\big(h_s(x),s\big),
\]
such that $h_0=\mathrm{id}$, each $h_s$ preserves the strata of $\Scal$, and \[h_s\big((V_0\cap M)\setminus N\big)=(V_s\cap M)\setminus N.\]
\end{theo}

The set $N$ meets every $V_s$ only along the set $V_0\cap\varphi^{-1}(0)$ of the family (see \cref{rem_baselocus}), so $h_s$ restricts to a homeomorphism from $(V_0\cap M)\setminus\varphi^{-1}(0)$ onto $(V_s\cap M)\setminus\varphi^{-1}(0)$. However, $h_s$ cannot, in general, be extended across $N$ by soft arguments (see \cref{rem_noextension}): to obtain an embedded trivialization, one needs to control $G$ over $\Delta_G$, and this is provided by Thom's $a_G$ condition. Following L\^e, Nu\~no-Ballesteros and Seade \cite[Definition~5.6 and Theorem~5.7]{LeNunoSeade2025}, a convenient hypothesis guaranteeing this condition is that the restriction of $\varphi$ to $V_0$ has an isolated critical point, that is,
\begin{equation}\label{eq:H}
\Sigma\big(\varphi|_{V_0}\big)=\{0\},
\end{equation}
where $\Sigma(\varphi|_{V_0})$ denotes the germ at $0$ of the stratified critical locus of $\varphi|_{V_0}$ (see \cref{def:H}). Our main result is the following.

\begin{theo}\label{theo_main}
Let $f_s:(X,0)\to(\C,0)$ be a separable family of analytic function-germs, with $s\in\C$.
If $(f_s)$ has a uniform critical set $\Sigma=\{0\}$ with respect to some Whitney stratification $\Scal$ of $X$ and $\varphi|_{V_0}$ has an isolated critical point, i.e., \eqref{eq:H} holds, then the family is embeddedly topologically trivial: there exist a neighborhood $M$ of $0$ in $X$, a neighborhood $D$ of $0$ in $\C$ and a homeomorphism
\[
M\times D\to M\times D,\qquad (x,s)\mapsto\big(h_s(x),s\big),
\]
such that $h_0=\mathrm{id}$, $h_s(0)=0$ and $h_s(V_0\cap M)=V_s\cap M$. In particular, the family has constant embedded topological type: for every $s\in D$, $h_s$ defines a homeomorphism of germs $(X,0)\to(X,0)$ sending $(V_0,0)$ onto $(V_s,0)$.
\end{theo}

Condition \eqref{eq:H} only involves the critical points of $\varphi$ on $V_0$, and $\varphi$ itself may have a non-isolated critical locus (see \cref{ex_1}). Moreover, since $f_0$ has an isolated critical point, \eqref{eq:H} implies that $G$ is transverse to small spheres off $\Delta_G$ (\cref{prop_Thom}), so \cref{theo_A} applies as well. On the other hand, \eqref{eq:H} is not necessary for the conclusion of \cref{theo_main} to hold (see \cref{ex_2}).

As an immediate application, we deduce that separable $\mu$-constant deformations of isolated complete intersection singularities (ICIS) are embeddedly trivial, provided $f_0$ and $\varphi$ define an ICIS of codimension two on $X$, which is equivalent to \eqref{eq:H} in this setting (see \cref{rem_H}).

\begin{cor}\label{cor_ICIS}
Let $(X,0)\subset(\C^N,0)$ be the germ of an ICIS and let $f_s:(X,0)\to(\C,0)$ be a separable family of isolated singularity function-germs such that the germ $X\cap f_0^{-1}(0)\cap\varphi^{-1}(0)$, with the structure defined by the equations of $X$, $f_0$ and $\varphi$, is an ICIS of dimension $\dim X-2$.
If $\mu(f_s)$ is constant for $|s|$ small, then $(f_s)$ is embeddedly topologically trivial; in particular, it has constant embedded topological type.
\end{cor}

The proofs build on the concepts of \textit{$\Delta$-regularity} and \textit{uniform critical set} from \cite{Menegon2023}, which we briefly recall in \cref{section_2}. \cref{theo_A} is proved in \cref{section_3} by lifting, through $G$, vector fields defined in the complement of $\Delta_G$. \cref{theo_main} is proved in \cref{section_4}, where we use \eqref{eq:H} to show that $G$ is a Thom map, which allows us to lift the same vector fields also over $\Delta_G$. In both sections we assume that the family is linear, and in \cref{section_5} we deduce the separable case by the change of parameter $t=\rho(s)$.

\section{\texorpdfstring{$\Delta$}{Delta}-regularity and uniform critical set}
\label{section_2}

We now recall two conditions, introduced in \cite{Menegon2023}, that provide control over the singular behavior of families of function-germs:
\emph{$\Delta$-regularity} and the existence of a \emph{uniform critical set}.

Let $(X,0)\subset(\CC^N,0)$ be a germ of complex analytic space and consider a separable family
\[
F(z,s)= f_0(z)+\rho(s)\,\phi(z),
\]
where $\rho:(\CC,0)\to(\CC,0)$ and $\phi:(\CC^N,0)\to(\CC,0)$ are analytic germs. We assume that $\rho$ is not identically zero, since otherwise the family is constant.
Restricting to $X$, we set
\[
\varphi:=\phi|_X:(X,0)\to(\CC,0), \qquad
f_s(x)= f_0(x)+\rho(s)\,\varphi(x)\ \ (x\in X).
\]

Since $(X,0)$ may have a complicated singular locus, fix a Whitney (equivalently, $(w)$-regular) stratification $\Scal=(S_\alpha)_{\alpha\in\Lambda}$ of $X$, and assume, as usual, that $\{0\}$ is a stratum.
For an analytic map $A$ defined on $X$, we denote by $\Sigma(A)$ its set of stratified critical points with respect to $\Scal$, i.e., the set of points $x\in X$ such that the restriction of $A$ to the stratum containing $x$ is not a submersion at $x$. In particular, $\Sigma(f_s)\subset X$ is the set of stratified critical points of $f_s$.

\medskip
Following \cite{Menegon2023}, define the auxiliary map
\[
G:(X,0)\to(\CC^2,0), \qquad z\mapsto \big(f_0(z),\varphi(z)\big),
\]
and its discriminant
\[
\Delta_G := G\big(\Sigma(G)\big),
\]
computed on a small representative of $G$.
We also consider
\[
\Gcal:(X\times \C,0)\to(\C^3,0), \qquad (z,s)\mapsto \big(f_0(z),\varphi(z),s\big).
\]

For $t\in\CC$, set
\[
H_t:=\{(y_0,y_1)\in\CC^2 \,;\, y_0+t\,y_1=0\}.
\]
Then
\begin{equation}\label{eq:preimG}
G^{-1}(H_{\rho(s)}) \;=\; f_s^{-1}(0) \;=\; V_s.
\end{equation}

\begin{definition}\label{def:deltareguniform}
The separable family $F$ is \emph{$\Delta$-regular} (with respect to $\Scal$) if the line $H_0$ is not a limit of secant lines of $\Delta_G$ at $(0,0)\in\CC^2$.

The family $F$ has a \emph{uniform critical set} (with respect to $\Scal$) if there exist $\varepsilon,\delta>0$ such that
\[
\Sigma(f_s)\cap \Bbb_\varepsilon \;=\; \Sigma(f_0)\cap \Bbb_\varepsilon
\qquad\text{for all } s\in \Dbb_\delta.
\]
In this case we set $\Sigma:=\Sigma(f_0)\cap \Bbb_\varepsilon$ and call $\Sigma$ the \emph{uniform critical set} of $F$.
\end{definition}

\begin{remark}\label{remark_delta}
Equivalently, $F$ is $\Delta$-regular if there exist an open disc $U\subset\CC^2$ centered at $0$, that is, $U=\{y\in\C^2\,;\,\norm{y}<c\}$ for some $c>0$, and $\eta>0$ such that
\[
H_t\cap \Delta_G \cap U \subset \{(0,0)\}
\qquad \text{whenever }|t|<\eta.
\]
Since $H_t=\{y_0=-t\,y_1\}$, this means that $\Delta_G\cap U$ is contained in the cone $\{|y_0|\geq\eta|y_1|\}$, and shrinking $U$ or the representative of $G$ only makes $\Delta_G\cap U$ smaller, so this property persists.
In terms of the unfolding of a linear family, that is, if $\rho(s)=s$, the hypersurface $\Hcal:=\{(y,s)\in\C^2\times\C\,;\,y\in H_s\}$ satisfies
\[
\Hcal\cap \Delta_\Gcal \cap (U\times \Dbb_{\eta}) \subseteq \{(0,0)\}\times\CC,
\]
where $\Delta_\Gcal=\Delta_G\times \C$ (see \cref{fig:Family}).
\end{remark}

The next proposition relates these two notions. It was proved in \cite{Menegon2023} for linear families, and the separable case follows from the linear one (see \cref{section_5}).

\begin{prop}\label{prop_MJA}
If a family $F$ has a uniform critical set (with respect to $\Scal$), then it is $\Delta$-regular (with respect to $\Scal$).
\end{prop}

The converse of \cref{prop_MJA} does not hold.

\begin{example}[cf.\ {\cite[Example~2.5]{Menegon2023}}]\label{ex_converse}
Let $X=\C^2$ with the stratification $\{\C^2\setminus\{0\},\{0\}\}$, $f_0=xy$ and $\varphi=x$, and consider the linear family $f_s=f_0+s\varphi=x(y+s)$.

The Jacobian matrix of $G=(f_0,\varphi)$ is
\[
\begin{pmatrix} y & x\\ 1 & 0\end{pmatrix},
\]
whose determinant is $-x$. Hence $\Sigma(G)=\{x=0\}$, on which both $f_0$ and $\varphi$ vanish, that is, $\Sigma(G)\subset G^{-1}(0)$. Therefore $\Delta_G=\{0\}$, so $H_t\cap\Delta_G=\{0\}$ for every $t$ and the family is $\Delta$-regular.

On the other hand, the only critical point of $f_s=x(y+s)$ in $\C^2$ is $p_s=(0,-s)$. Since $\{0\}$ is a stratum, $0\in\Sigma(f_s)$ for every $s$, and therefore, for $0<|s|<\e$,
\[
\Sigma(f_s)\cap\Bbb_\e=\{0,p_s\}\neq\{0\}=\Sigma(f_0)\cap\Bbb_\e.
\]
Thus the family does not have a uniform critical set. 

Note that $V_s=\{x=0\}\cup\{y=-s\}$ is smooth at $0$ for $s\neq0$, whereas $V_0$ is a node. Hence $\Delta$-regularity alone does not imply constant topological type.
\end{example}

\begin{figure}
\centering
\definecolor{figblue}{rgb}{0,0.678,0.937}%
\definecolor{figred}{rgb}{0.929,0.110,0.141}%
\definecolor{figgray}{rgb}{0.655,0.649,0.650}%
\definecolor{figblack}{rgb}{0.137,0.122,0.125}%
\begin{tikzpicture}[x=0.8pt,y=0.8pt,font=\small,
  axis/.style={figgray,line width=0.4pt},
  line/.style={figblue,line width=1.2pt},
  branch/.style={figred,line width=1.2pt},
  cone/.style={figblack!55,densely dashed,line width=0.5pt}]
\path[use as bounding box] (0,0) rectangle (439.37,255.12);
\coordinate (O1) at (69.11,184.76);
\coordinate (O2) at (218.08,126.63);
\coordinate (O3) at (357.59,69.45);
\coordinate (A1) at (86.14,248.07);  \coordinate (B1) at (52.23,121.49);
\coordinate (A2) at (218.39,187.38); \coordinate (B2) at (217.91,61.07);
\coordinate (A3) at (340.70,132.76); \coordinate (B3) at (374.62,6.18);
\foreach \i in {1,2,3}{
  \path (O\i) ++(64.50,36.86) coordinate (U\i);
  \path (O\i) ++(64.50,-36.86) coordinate (D\i);
}
\fill[figred,fill opacity=0.2]
  (O1) .. controls +(23.36,0) and +(-23.54,-23.54) .. (U1) -- (U2)
  .. controls +(-23.54,-23.54) and +(23.36,0) .. (O2)
  .. controls +(23.36,0) and +(-23.54,23.54) .. (D2) -- (D1)
  .. controls +(-23.54,23.54) and +(23.36,0) .. (O1) -- cycle;
\fill[figblue,fill opacity=0.2] (B1) -- (A1) -- (A2) -- (B2) -- cycle;
\fill[figred,fill opacity=0.2]
  (O2) .. controls +(23.36,0) and +(-23.54,-23.54) .. (U2) -- (U3)
  .. controls +(-23.54,-23.54) and +(23.36,0) .. (O3)
  .. controls +(23.36,0) and +(-23.54,23.54) .. (D3) -- (D2)
  .. controls +(-23.54,23.54) and +(23.36,0) .. (O2) -- cycle;
\fill[figblue,fill opacity=0.2] (B2) -- (A2) -- (A3) -- (B3) -- cycle;
\draw[axis] (7.45,185.25) -- (134.71,184.78);
\draw[axis] (69.19,119.26) -- (69.19,250.31);
\draw[axis] (152.39,126.59) -- (283.43,126.59);
\draw[axis] (292.14,69.47) -- (423.19,69.47);
\draw[axis] (357.66,3.95) -- (357.66,134.99);
\foreach \i in {1,2,3}{
  \draw[cone] (O\i) ++(-28.6,-63.5) -- ++(57.2,127);
  \draw[cone] (O\i) ++(28.6,-63.5) -- ++(-57.2,127);
}
\draw[line] (B1) -- (A1);
\draw[line] (B2) -- (A2);
\draw[line] (B3) -- (A3);
\foreach \i in {1,2,3}{
  \draw[branch] (O\i) .. controls +(23.36,0) and +(-23.54,-23.54) .. (U\i);
  \draw[branch] (O\i) .. controls +(23.36,0) and +(-23.54,23.54) .. (D\i);
}
\draw[figblack,line width=1.2pt] (O1) -- (O2) -- (O3);
\node[figblue] at (107.46,82.54) {$\Hcal\cap\Ucal_\Theta$};
\node[figblue] at (209.62,52.31) {$H_0$};
\node[figblue] at (389.50,15.69) {$H_s$};
\node[figred] at (345.00,158.47) {$\Delta_\Gcal\cap\Ucal_\Theta$};
\node[figred] at (423.19,87.23) {$\Delta_G$};
\node[figblack!75] at (160.38,140.28) {$\Theta$};
\end{tikzpicture}%
\caption{Depiction of the different objects in our proof, over a segment $\Theta$ of parameters through $0$. Here $\Ucal_\Theta=U\times\Theta$, and the dashed lines bound, in each plane $\C^2\times\{s\}$, the cone $\{|y_0|\leq r|y_1|\}$ of Step~1 in \cref{section_3}, which contains the lines $H_s$ and meets $\Delta_G$ only at the origin.}
\label{fig:Family}
\end{figure}

\section{Trivialization outside the discriminant}
\label{section_3}

In this section we prove \cref{theo_A} for linear families; the separable case is deduced in \cref{section_5}. Thus, in this section and the next one, we assume that $\rho(s)=s$, that is,
\[
f_s=f_0+s\,\varphi\qquad\text{and}\qquad V_s=G^{-1}(H_s).
\]
We also assume that $(f_s)$ is $\Delta$-regular with respect to $\Scal$. Let $U$ and $\eta$ be as in \cref{remark_delta}. We fix $\e>0$ small enough so that $\mathbb{S}_\e$ is transverse to every stratum of $\Scal$, which is possible by the curve selection lemma \cite[Lemma~3.1]{Milnor1968}, and consider the representative $G:X\cap\overline{\Bbb}_\e\to\C^2$, its discriminant $\Delta_G$ and
\[
N:=G^{-1}(\Delta_G).
\]
By \cref{eq:preimG}, the fibers $V_s$ are described as preimages under $G$. Outside $N$ the map $G$ is a stratified submersion. To obtain a trivialization, we also need the fibers of $G$ to be transverse to the boundary of the ball.

\begin{definition}\label{def:ML}
We say that $G$ is \emph{transverse to small spheres off $\Delta_G$} if, for every sufficiently small $\e>0$, there exists $\nu>0$ such that $\mathbb{S}_\e$ is transverse to $G^{-1}(y)\cap S_\alpha$ for every $y\in\Dbb_\nu\setminus\Delta_G$ and every stratum $S_\alpha$ of $\Scal$, where $\Dbb_\nu\subset\CC^2$ is the open disc of radius $\nu$ centered at $0$.
\end{definition}

In this case, $G$ restricts to a proper stratified submersion from $X\cap\overline{\Bbb}_\e\cap G^{-1}(\Dbb_\nu\setminus\Delta_G)$ onto $\Dbb_\nu\setminus\Delta_G$ (where we also consider the strata $S_\alpha\cap\mathbb{S}_\e$), hence $G$ induces a locally trivial fibration by Thom--Mather's first isotopy lemma \cite[Ch.~II]{Gibson1976a} (cf. \cite{Mather2012}). This is the Milnor--L\^e fibration of $G$ outside its discriminant, see \cite[Definition~5.1]{LeNunoSeade2025}.

\begin{remark}\label{rem_ML}
By (the proof of) \cite[Lemma~5.2]{LeNunoSeade2025} (see also the proof of \cref{prop_Thom}), the condition of \cref{def:ML} holds whenever $G$ satisfies Thom's $a_G$ condition with respect to a Whitney stratification of $X$ refining $\Scal$ in which $G^{-1}(0)$ is a union of strata, and the stratification $\big\{U\setminus\Delta_G,\Delta_G\setminus\{0\},\{0\}\big\}$ of the target. In particular, it holds under \eqref{eq:H} if $f_0$ has an isolated critical point (\cref{prop_Thom}). 

However, it does not follow from $\Delta$-regularity nor from the existence of a uniform critical set as the following example shows. 
\end{remark}

\begin{example}\label{ex_LeSabbah}
We consider L\^e--Sabbah classical example of a map without Thom's condition \cite{Sabbah1983}. Let $X=\C^3$ with the stratification $\{\C^3\setminus\{0\},\{0\}\}$, $f_0=y$ and $\varphi=x^2-y^2z$, $G=(f_0,\varphi)$. Since $\partial f_s/\partial y=1-2syz\neq0$ near $0$, the family has a uniform critical set $\Sigma=\{0\}$. Moreover, $\Sigma(G)$ is the $z$-axis, so $\Delta_G=\{0\}$. 

However, at every point $p=(0,a,b)$ with $a,b\neq0$ small, the kernel of $d_pG$ is the line $\C\times\{0\}\times\{0\}$. Hence the fiber $G^{-1}\big(G(p)\big)$ is a curve tangent to the sphere that contains $p$. Since we can take $p\in\mathbb{S}_\e$ so that the values $G(p)=(a,-a^2b)$ are not in $\Delta_G$, the condition of \cref{def:ML} fails for every $\e$.
\end{example}

We proceed now to prove \cref{theo_A}, which we separate by steps. To construct a trivialization, we first need good control in the target space.

\medskip
\noindent
\textbf{Step 1. The parameter.}
We write $s=s_1+i\,s_2$, identifying $\C$ with $\RR^2$, and set $e_1:=1$ and $e_2:=i$. Fix $0<r<\eta$. Hence, for $s\in\Dbb_r$, the line $H_s$ is contained in the closed cone $\{|y_0|\leq r|y_1|\}$, which meets $\Delta_G\cap U$ only at $0$ by \cref{remark_delta} (see \cref{fig:Family}). We will move the lines $H_s$, $s\in\Dbb_r$, onto the line $H_0$ while keeping $\Delta_G$ fixed.

\medskip
\noindent
\textbf{Step 2. Vector fields adapted to $\Hcal$.}
By \cref{remark_delta}, inside $U\times\Dbb_{r}$ the intersection $\Hcal\cap\Delta_\Gcal$ lies entirely in the axis $\{(0,0)\}\times\Dbb_{r}$ (see \cref{fig:Family}). We want two vector fields $\vec v_1,\vec v_2$ on $U\times\Dbb_r$, with $s$-components $\partial/\partial s_1$ and $\partial/\partial s_2$, respectively, which are tangent to both $\Delta_\Gcal$ and $\Hcal\cap(U\times\Dbb_r)$ (see \cref{fig:vfields,fig:vslice}). We will use them to reach any parameter $s=s_1+is_2$ from the origin moving first to $s_1$. We obtain these vector fields in two stages. First, \cref{prop_vec_v} gives a convenient map $\Psi$ and, then, we define $\vec v_j:=d\Psi^{-1}(\partial/\partial s_j)$.

\begin{proposition}\label{prop_vec_v}
There is a homeomorphism
\[
\Psi(y,s)=\big(\psi_s(y),s\big)
\]
of $U\times\Dbb_r$, which restricts to a diffeomorphism of $(U\setminus\{0\})\times\Dbb_r$, such that:
\begin{enumerate}[(i)]
\item $\psi_0=\mathrm{id}$, $\norm{\psi_s(y)}=\norm{y}$ and $\psi_s(\lambda y)=\lambda\,\psi_s(y)$ for $y\in U$ and $\lambda\in(0,1]$;
\item $\psi_s(H_s\cap U)=H_0\cap U$ for every $s\in\Dbb_r$, that is, $\Psi$ sends $\Hcal\cap(U\times\Dbb_{r})$ to $(H_0\cap U)\times\Dbb_{r}$;
\item $\Psi$ fixes $\Omega\times\Dbb_r$, where $\Omega\subset U$ is an open neighborhood of $(\Delta_G\cap U)\setminus\{0\}$.
\end{enumerate}
\end{proposition}

\begin{proof}
We construct $\Psi$ on $\mathbb{S}^3\times\Dbb_r$, where $\mathbb{S}^3\subset\C^2$ is the unit sphere, and extend it by homogeneity. Set
\[
W:=\{u\in\mathbb{S}^3\,;\,|u_0|<r\,|u_1|\},\quad\text{so that}\quad\overline W=\{u\in\mathbb{S}^3\,;\,|u_0|\leq r\,|u_1|\}.
\]
Since $|s|<r$ for $s\in\Dbb_r$, the circles $K_s:=H_s\cap\mathbb{S}^3$, $s\in\Dbb_r$, are contained in $W$, and
\[
\mathcal K:=\{(u,s)\in\mathbb{S}^3\times\Dbb_r\,;\,u\in K_s\}
\]
is a closed submanifold of $\mathbb{S}^3\times\Dbb_r$.

Indeed, since $H_s$ is spanned by the unit vector $w_s:=(-s,1)/\norm{(-s,1)}$, the map
\[
\mathbb{S}^1\times\Dbb_r\to\mathcal K,\qquad(\lambda,s)\mapsto(\lambda\,w_s,s),
\]
is a diffeomorphism which commutes with the projections onto $\Dbb_r$. This diffeomorphism sends the vector field $\partial/\partial s_j$ on $\mathbb{S}^1\times\Dbb_r$ to a vector field on $\mathcal K$ of the form $(\xi_j,\partial/\partial s_j)$, with $\xi_j(u,s)\in T_u\mathbb{S}^3$, $j=1,2$. We can now extend $\xi_j$ smoothly to $\mathbb{S}^3\times\Dbb_r$ (similarly to the proof of the isotopy extension theorem \cite[Ch.~8, Theorem~1.3]{Hirsch1976}). Then, we multiply it by a smooth function of $u$ which is equal to $1$ on $\overline W$ and vanishes outside a compact set $L\subset\{u\in\mathbb{S}^3\,;\,|u_0|<\eta\,|u_1|\}$, which exists because $r<\eta$ (in \cref{fig:vslice}, the cone over $L$ is bounded by the dotted lines). 

We keep the notation for the resulting vector field on $\mathbb{S}^3\times\Dbb_r$, that is, $(\xi_j,\partial/\partial s_j)$. Note that it is still tangent to $\mathcal K$. Since $\mathbb{S}^3$ is compact and $\Dbb_r$ is convex, its flow is defined at $(u,s)$ for every time $\tau\in\RR$ with $s+\tau e_j\in\Dbb_r$, and it has the form $(u,s)\mapsto\big(g^j_{\tau,s}(u),s+\tau e_j\big)$, where $g^j_{\tau,s}$ is a diffeomorphism of $\mathbb{S}^3$ which sends $K_s$ onto $K_{s+\tau e_j}$ and is the identity outside $L$. 

If $s=s_1+i\,s_2\in\Dbb_r$, then also $s_1\in\Dbb_r$, so
\[
g_s:=g^2_{s_2,s_1}\circ g^1_{s_1,0}
\]
is defined and it is a diffeomorphism of $\mathbb{S}^3$ which, together with its inverse $g_s^{-1}=g^1_{-s_1,s_1}\circ g^2_{-s_2,s}$, depends smoothly on $s$, and satisfies $g_0=\mathrm{id}$, $g_s(K_0)=K_s$ and $g_s=\mathrm{id}$ outside $L$. We use $g_s$ to define $\psi_s\colon U\to U$ by
\[
\psi_s(y):=\norm{y}\,g_s^{-1}\big(y/\norm{y}\big)\ \ (y\neq0),\qquad\psi_s(0):=0.
\]
Then $\Psi(y,s)=\big(\psi_s(y),s\big)$ is a diffeomorphism of $(U\setminus\{0\})\times\Dbb_r$, with inverse $(y,s)\mapsto\big(\norm{y}\,g_s(y/\norm{y}),s\big)$, and it is a homeomorphism of $U\times\Dbb_r$, since $\norm{\psi_s(y)}=\norm{y}$. Moreover, (i) holds by construction, since $g_0=\mathrm{id}$. Since $H_s\cap U$ is the cone over $K_s=H_s\cap\mathbb{S}^3$, (ii) follows from $g_s(K_0)=K_s$. Finally, (iii) holds for $\Omega:=\{y\in U\setminus\{0\}\,;\,y/\norm{y}\notin L\}$, which is open and contains $(\Delta_G\cap U)\setminus\{0\}$, because $\Delta_G\cap U\subset\{|y_0|\geq\eta|y_1|\}$ by \cref{remark_delta}.
\end{proof}

Then, we consider the vector fields
\[
\vec{v}_j\coloneqq d\Psi^{-1}\left(\frac{\partial}{\partial s_j}\right),\qquad j=1,2,
\]
which are as required by construction of $\Psi$. More precisely, set
\[\theta^j_\tau(y,s):=(y,s+\tau e_j)\]
for $\tau\in\RR$ with $s+\tau e_j\in\Dbb_r$. Then $\vec v_j$ is the vector field whose flow is $\Psi^{-1}\circ\theta^j_\tau\circ\Psi$, defined at $(y,s)$ whenever $s+\tau e_j\in\Dbb_r$, so it is smooth on $(U\setminus\{0\})\times\Dbb_r$, and it is tangent to $\Hcal\cap(U\times\Dbb_{r})$ by \cref{prop_vec_v}\,(ii). Moreover, $\vec v_j=\partial/\partial s_j$ on $\{0\}\times\Dbb_r$ and, by \cref{prop_vec_v}\,(iii), on $\Omega\times\Dbb_r$ (see \cref{fig:vfields,fig:vslice}). Finally, $\vec v_j$ is continuous, since its component in $\C^2$ is positively homogeneous of degree one in $y$ by \cref{prop_vec_v}\,(i).

\begin{figure}
\centering
\definecolor{figblue}{rgb}{0,0.678,0.937}%
\definecolor{figred}{rgb}{0.929,0.110,0.141}%
\definecolor{figgray}{rgb}{0.655,0.649,0.650}%
\definecolor{figblack}{rgb}{0.137,0.122,0.125}%
\tikzset{
  axis/.style={figgray,line width=0.4pt},
  line/.style={figblue,line width=1.2pt},
  branch/.style={figred,line width=1.2pt},
  cone/.style={figblack!55,densely dashed,line width=0.5pt},
  flowH/.style={figblue!60!black,line width=0.7pt,
    postaction={decorate},decoration={markings,
      mark=at position 0.27 with {\arrow{Stealth[length=4pt]}},
      mark=at position 0.77 with {\arrow{Stealth[length=4pt]}}}},
  flowD/.style={figred!65!black,line width=0.7pt,dash pattern=on 3pt off 1.2pt,
    postaction={decorate},decoration={markings,
      mark=at position 0.27 with {\arrow{Stealth[length=4pt]}},
      mark=at position 0.77 with {\arrow{Stealth[length=4pt]}}}},
}
\begin{tikzpicture}[x=0.8pt,y=0.8pt,font=\small]
\path[use as bounding box] (0,-12) rectangle (439.37,255.12);
\coordinate (O1) at (69.11,184.76);
\coordinate (O2) at (218.08,126.63);
\coordinate (O3) at (357.59,69.45);
\coordinate (A1) at (86.14,248.07);  \coordinate (B1) at (52.23,121.49);
\coordinate (A2) at (218.39,187.38); \coordinate (B2) at (217.91,61.07);
\coordinate (A3) at (340.70,132.76); \coordinate (B3) at (374.62,6.18);
\foreach \i in {1,2,3}{
  \path (O\i) ++(64.50,36.86) coordinate (U\i);
  \path (O\i) ++(64.50,-36.86) coordinate (D\i);
  \path (O\i) ++(47.78,21.17) coordinate (R\i);
  \path (O\i) ++(22.91,4.76) coordinate (RR\i);
  \path (O\i) ++(59.64,-32.07) coordinate (S\i);
  \path ($(O\i)!0.45!(A\i)$) coordinate (P\i);
  \path ($(O\i)!0.85!(A\i)$) coordinate (PP\i);
  \path ($(O\i)!0.45!(B\i)$) coordinate (Q\i);
  \path ($(O\i)!0.85!(B\i)$) coordinate (QQ\i);
}
\fill[figred,fill opacity=0.13]
  (O1) .. controls +(23.36,0) and +(-23.54,-23.54) .. (U1) -- (U2)
  .. controls +(-23.54,-23.54) and +(23.36,0) .. (O2)
  .. controls +(23.36,0) and +(-23.54,23.54) .. (D2) -- (D1)
  .. controls +(-23.54,23.54) and +(23.36,0) .. (O1) -- cycle;
\fill[figblue,fill opacity=0.13] (B1) -- (A1) -- (A2) -- (B2) -- cycle;
\fill[figred,fill opacity=0.13]
  (O2) .. controls +(23.36,0) and +(-23.54,-23.54) .. (U2) -- (U3)
  .. controls +(-23.54,-23.54) and +(23.36,0) .. (O3)
  .. controls +(23.36,0) and +(-23.54,23.54) .. (D3) -- (D2)
  .. controls +(-23.54,23.54) and +(23.36,0) .. (O2) -- cycle;
\fill[figblue,fill opacity=0.13] (B2) -- (A2) -- (A3) -- (B3) -- cycle;
\draw[axis] (7.45,185.25) -- (134.71,184.78);
\draw[axis] (69.19,119.26) -- (69.19,250.31);
\draw[axis] (152.39,126.59) -- (283.43,126.59);
\draw[axis] (292.14,69.47) -- (423.19,69.47);
\draw[axis] (357.66,3.95) -- (357.66,134.99);
\draw[line] (B1) -- (A1);
\draw[line] (B2) -- (A2);
\draw[line] (B3) -- (A3);
\foreach \i in {1,2,3}{
  \draw[branch] (O\i) .. controls +(23.36,0) and +(-23.54,-23.54) .. (U\i);
  \draw[branch] (O\i) .. controls +(23.36,0) and +(-23.54,23.54) .. (D\i);
}
\foreach \X in {P,PP,Q,QQ}{
  \draw[flowH] plot[smooth,tension=0.6] coordinates {(\X1) (\X2) (\X3)};
}
\foreach \X in {R,RR,S}{
  \draw[flowD] plot[smooth,tension=0.6] coordinates {(\X1) (\X2) (\X3)};
}
\draw[figblack,line width=1.2pt,postaction={decorate},decoration={markings,
  mark=at position 0.27 with {\arrow{Stealth[length=5pt]}},
  mark=at position 0.77 with {\arrow{Stealth[length=5pt]}}}] (O1) -- (O2) -- (O3);
\node[figblue] at (209.62,52.31) {$H_0$};
\node[figblue] at (389.50,15.69) {$H_s$};
\node[figblue] at (40,108) {$H_{-s}$};
\node[figred] at (423.19,87.23) {$\Delta_G$};
\node[figblue!60!black] at ($(QQ2)+(48,-26)$) {$\vec v_1$};
\node[figred!65!black,anchor=south west] at (342,150) {$\vec v_1=\partial/\partial s_1$};
\draw[figred!65!black,line width=0.4pt] (379,150) -- ($(R2)!0.82!(R3)$);
\node[figblack!75,fill=white,fill opacity=0.85,text opacity=1,inner sep=1pt,rounded corners=2pt] at (160.38,140.28) {$\Theta$};
\end{tikzpicture}
\caption{Integral curves of the vector field $\vec v_1$ of Step~2 in \cref{section_3} over a real segment $\Theta$ of parameters through $0$. On $\Hcal$ (solid) they join points of $H_{-s}$, $H_0$ and $H_s$ at the same distance from the axis, whereas on $\Delta_\Gcal$ (dashed) and on the axis $\{0\}\times\Theta$ they are parallel to $\Theta$, since $\vec v_1=\partial/\partial s_1$ there. The vector field $\vec v_2$ behaves in the same way for the imaginary direction of the parameter.}
\label{fig:vfields}
\end{figure}

\begin{figure}
\centering
\definecolor{figblue}{rgb}{0,0.678,0.937}%
\definecolor{figred}{rgb}{0.929,0.110,0.141}%
\definecolor{figgray}{rgb}{0.655,0.649,0.650}%
\definecolor{figblack}{rgb}{0.137,0.122,0.125}%
\tikzset{
  axis/.style={figgray,line width=0.4pt},
  line/.style={figblue,line width=1.2pt},
  branch/.style={figred,line width=1.2pt},
  cone/.style={figblack!55,densely dashed,line width=0.5pt},
}
\begin{tikzpicture}[x=1pt,y=1pt,font=\small]
\fill[figgray!25] (0,0) -- (80,80) -- (80,-80) -- cycle;
\fill[figgray!25] (0,0) -- (-80,80) -- (-80,-80) -- cycle;
\draw[axis] (-85,0) -- (85,0) node[right,figblack] {$y_0$};
\draw[axis] (0,-85) -- (0,85) node[above,figblack] {$y_1$};
\draw[figblack!45,dotted,line width=0.6pt] (-80,-80) -- (80,80);
\draw[figblack!45,dotted,line width=0.6pt] (80,-80) -- (-80,80);
\draw[cone] (-36,-80) -- (36,80);
\draw[cone] (36,-80) -- (-36,80);
\draw[branch] (0,0) .. controls (23.36,0) and (40.96,13.32) .. (64.5,36.86);
\draw[branch] (0,0) .. controls (23.36,0) and (40.96,-13.32) .. (64.5,-36.86);
\node[figred] at (75,44) {$\Delta_G$};
\draw[line] (-21.4,80) -- (21.4,-80);
\node[figblue,anchor=north] at (21.4,-80) {$H_s$};
\foreach \rr in {22,38,54,70}{
  \foreach \pp in {-40,-32,...,40}{
    \pgfmathsetmacro{\ang}{abs(\pp)}
    \pgfmathsetmacro{\cut}{(\ang<=24.2) ? 1 : ((\ang>=45) ? 0 : (45-\ang)/(45-24.2))}
    \pgfmathsetmacro{\len}{0.32*\cut*\rr}
    \ifdim\len pt>0.5pt
      \draw[-{Stealth[length=3.2pt]},figblue!60!black,line width=0.5pt]
        ({-\rr*sin(\pp)},{\rr*cos(\pp)}) -- ++({-\len*cos(\pp)},{-\len*sin(\pp)});
      \draw[-{Stealth[length=3.2pt]},figblue!60!black,line width=0.5pt]
        ({\rr*sin(\pp)},{-\rr*cos(\pp)}) -- ++({\len*cos(\pp)},{\len*sin(\pp)});
    \fi
  }
}
\fill (0,0) circle (1.2pt);
\end{tikzpicture}
\caption{The component in $\C^2$ of the vector field $\vec v_1$ of Step~2 in \cref{section_3} on a plane $\C^2\times\{s\}$ with $s>0$, for one choice of the extension of $\xi_1$ in the proof of \cref{prop_vec_v}. It moves the line $H_s$ towards $H_{s+\tau}$, $\tau>0$, inside the cone $\{|y_0|\leq r|y_1|\}$ (dashed), it vanishes outside the cone over $L$ (dotted), that is, on $\Omega$ (shaded), which contains $\Delta_G\setminus\{0\}$, and it grows linearly with $\|y\|$.}
\label{fig:vslice}
\end{figure}

\medskip
\noindent
\textbf{Step 3. Lifting the vector fields.}
Assume now that $G$ is transverse to small spheres off $\Delta_G$, with $\e$ small enough for \cref{def:ML} to apply, and shrink $U$ if necessary to be contained in the disc $\Dbb_\nu$ of \cref{def:ML} ($\Psi$ restricts to it by \cref{prop_vec_v}\,(i)). Set
\[
M:=X\cap\overline{\Bbb}_\e\cap G^{-1}(U),
\]
which is a neighborhood of $0$ in $X$, and stratify $M\setminus N$ by the sets $(S_\alpha\cap\mathring{\Bbb}_\e)\setminus N$ and $(S_\alpha\cap\mathbb{S}_\e)\setminus N$, $S_\alpha\in\Scal$; this is a Whitney stratification, since $\mathbb{S}_\e$ is transverse to the strata of $\Scal$ (see, e.g., \cite[Lemma~2.5]{GimenezConejero2023}). Then the restriction of the unfolding $\Gcal$,
\[
\Gcal:(M\setminus N)\times\Dbb_r\longrightarrow(U\setminus\Delta_G)\times\Dbb_r,\qquad (x,s)\longmapsto\big(G(x),s\big),
\]
is a stratified submersion: on the strata contained in $\mathring{\Bbb}_\e$ because $\Sigma(G)\subset N$, and on the strata contained in $\mathbb{S}_\e$ by \cref{def:ML}. Moreover, it is proper, since the preimage of a compact set $K\times L$ is the compact set $\big(X\cap\overline{\Bbb}_\e\cap G^{-1}(K)\big)\times L$. By \cref{prop_vec_v}\,(iii), the flows of $\vec v_1$ and $\vec v_2$ preserve the set $(U\setminus\Delta_G)\times\Dbb_r$, so these are smooth vector fields on $(U\setminus\Delta_G)\times\Dbb_r$ (recall that $0\in\Delta_G$).

Recall that a stratified submersion onto a manifold, with the trivial stratification, is a Thom map \cite[II.3.3]{Gibson1976a}, so we can use the lifting results \cite[Theorems~3.5, 3.7 and~3.9]{GimenezConejero2023}, which follow \cite[Ch.~II]{Gibson1976a}; or, alternatively, Verdier's lifting results \cite[Prop.~4.6]{Verdier1976}. 

Each vector field $\vec v_j$ on $(U\setminus\Delta_G)\times\Dbb_r$ admits a locally integrable stratified lift $\vec u_j$ to $(M\setminus N)\times\Dbb_r$, that is, $d\Gcal(\vec u_j)=\vec v_j$. Since $\Gcal$ is proper, the flow $\Phi^j_\tau$ of $\vec u_j$ is defined at $(x,s)$ whenever the flow of $\vec v_j$ is defined at $\Gcal(x,s)$, that is, whenever $s+\tau e_j\in\Dbb_r$ (cf.\ \cite[Theorem~3.10]{GimenezConejero2023}). 
Moreover, $\Phi^j_\tau$ is continuous, preserves the strata, and satisfies the $\Gcal$-relation between the flows
\begin{equation}\label{eq:commflows}
\Gcal\circ\Phi^j_\tau=\big(\Psi^{-1}\circ\theta^j_\tau\circ\Psi\big)\circ\Gcal.
\end{equation}

\medskip
\noindent
\textbf{Step 4. Conclusion.}

\begin{proof}[Proof of \cref{theo_A} for linear families]
Since the $s$-component of $\vec u_j$ is $\partial/\partial s_j$, the flow $\Phi^j_\tau$ maps $(M\setminus N)\times\{s\}$ onto $(M\setminus N)\times\{s+\tau e_j\}$ whenever $s,s+\tau e_j\in\Dbb_r$. Hence, for $s=s_1+i\,s_2\in\Dbb_r$, since also $s_1\in\Dbb_r$, we may define $h_s:M\setminus N\to M\setminus N$ by
\[
\big(h_s(x),s\big):=\Phi^2_{s_2}\big(\Phi^1_{s_1}(x,0)\big).
\]
It is clear that the map $(x,s)\mapsto\big(h_s(x),s\big)$ is a (stratified) homeomorphism. Moreover, $h_0=\mathrm{id}$, and $G\circ h_s=\psi_s^{-1}\circ G$ by \cref{prop_vec_v,eq:commflows}. Therefore, by \cref{prop_vec_v}\,(ii), for $s\in\Dbb_{r}$,
\begin{align*}
h_s\big((V_0\cap M)\setminus N\big)&=G^{-1}\big(\psi_s^{-1}(H_0\cap U)\big)\cap(M\setminus N)\\
&=G^{-1}\big(H_{s}\cap U\big)\cap(M\setminus N)=(V_s\cap M)\setminus N,
\end{align*}
which proves the theorem with $D=\Dbb_{r}$.
\end{proof}

\begin{remark}\label{rem_baselocus}
If $|t|<\eta$, then $H_t\cap\Delta_G\cap U\subset\{0\}$ by \cref{remark_delta}, so $G^{-1}(H_t)\cap N\cap M=G^{-1}(0)\cap M$. Hence, for $|s|<\eta$, $N$ meets $V_s$ only along the base locus $V_0\cap\varphi^{-1}(0)$ of the family, which is the same for every $s$, and the homeomorphism $h_s$ of \cref{theo_A} restricts to a homeomorphism
\[
(V_0\cap M)\setminus\varphi^{-1}(0)\longrightarrow(V_s\cap M)\setminus\varphi^{-1}(0).
\]
\end{remark}

\begin{remark}\label{rem_cone}
The hypothesis of \cref{def:ML} is only used to choose the lifts of $\vec v_1,\vec v_2$ tangent to $\mathbb{S}_\e$, so that their flows do not leave $M$. Since $\vec v_j=\partial/\partial s_j$ on $\Omega\times\Dbb_r$, over $G^{-1}(\Omega)$ one may take the trivial lift $(0,\partial/\partial s_j)$ and glue it to the previous one with a partition of unity. Then $h_s$ is the identity outside 
\[G^{-1}(U\setminus\Omega)\subset\bigcup_{|t|<\eta}G^{-1}(H_t).\]
 Hence, it suffices that $\mathbb{S}_\e$ be transverse to 
\[G^{-1}(H_t)\cap\varphi^{-1}(c)\cap S_\alpha\]
 for $|t|<\eta$ and $0<|c|$ small, that is, that the functions $\varphi|_{G^{-1}(H_t)}$ have Milnor fibrations of radius $\e$ uniformly in $t$. This weaker condition holds in L\^e--Sabbah's \cref{ex_LeSabbah}. We leave the long but straightforward computations of this fact to the reader.
\end{remark}

\begin{remark}\label{rem_noextension}
A homeomorphism defined on the complement of a subset of real codimension two does not extend, in general, to a homeomorphism of the whole space. For instance, $z\mapsto(1-|z|)\,z/|z|$ is a homeomorphism of the punctured disc $\{0<|z|<1\}$ that does not extend continuously to the origin. Therefore, extending $h_s$ across $N$ requires some control of the lifted vector fields near $N$. This control is provided by Thom's condition, as we show in the next section.
\end{remark}

\section{\texorpdfstring{Embedded triviality under a Thom condition}{Embedded triviality under a Thom condition}}
\label{section_4}

In this section we prove \cref{theo_main} and \cref{cor_ICIS} for linear families, that is, we keep assuming that $\rho(s)=s$ (see \cref{section_5} for the separable case). We assume that $(f_s)$ has a uniform critical set $\Sigma=\{0\}$ with respect to $\Scal$. Then $(f_s)$ is $\Delta$-regular by \cref{prop_MJA} and, since $f_0|_{S_\alpha}$ is a submersion on $S_\alpha\setminus\{0\}$, the germ $V_0\setminus\{0\}$ is stratified by the manifolds $V_0\cap S_\alpha$, $S_\alpha\in\Scal$.

\begin{definition}\label{def:H}
We say that $\varphi|_{V_0}$ has an \emph{isolated critical point} at $0$ if its stratified critical locus is $\Sigma(\varphi|_{V_0})=\{0\}$ as a germ, that is, if the restriction of $\varphi$ to each $V_0\cap S_\alpha$, $S_\alpha\neq\{0\}$, has no critical points near $0$. This is condition \eqref{eq:H}.
\end{definition}

Let $x\in V_0\cap S_\alpha$, $x\neq0$. Since $f_0|_{S_\alpha}$ is a submersion at $x$, the restriction $G|_{S_\alpha}$ is a submersion at $x$ if and only if $\varphi|_{V_0\cap S_\alpha}$ is a submersion at $x$. Therefore, condition \eqref{eq:H} is equivalent to
\begin{equation}\label{eq:H2}
\Sigma(G)\cap V_0=\{0\}.
\end{equation}

\begin{remark}\label{rem_H}

\begingroup
By $\Delta$-regularity, we have the equality of germs
\[
\Sigma(G)\cap V_0=\Sigma(G)\cap G^{-1}(0).
\]
Indeed, if $x\in\Sigma(G)\cap V_0$, then
\[
G(x)\in\Delta_G\cap H_0=\{0\},
\]
after choosing sufficiently small representatives.
The reverse inclusion follows from $G^{-1}(0)\subset V_0$.

Consequently, by \eqref{eq:H2}, condition \eqref{eq:H} is
equivalent to
\[
\Sigma(G)\cap G^{-1}(0)=\{0\}
\]
as germs.

Suppose now that $X\setminus\{0\}$ is smooth and
$\Scal=\{X\setminus\{0\},\{0\}\}$.
The preceding equivalence says that \eqref{eq:H} holds
if and only if $G|_{X\setminus\{0\}}$ is a submersion at
every point of $G^{-1}(0)\setminus\{0\}$ sufficiently close
to the origin.
In particular, when $X$ is an ICIS, the Jacobian criterion
shows that \eqref{eq:H} holds if and only if the equations
of $X$, together with $f_0$ and $\varphi$, define an ICIS
of the expected dimension $\dim X-2$, with the analytic
structure specified in \cref{cor_ICIS}.

The condition on the dimension cannot be dropped:
for $X=\C^3$, $f_0=z$ and $\varphi=xz$, the set
\[
\{f_0=\varphi=0\}=\{z=0\}
\]
is smooth, but it has dimension two rather than
$\dim X-2=1$. Moreover,
\[
\Sigma(G)=G^{-1}(0)=\{z=0\}=V_0,
\]
so \eqref{eq:H} does not hold, even though
$\Delta_G=\{0\}$ and the family is $\Delta$-regular.

\par
\endgroup

\end{remark}

The key point is that, under \eqref{eq:H}, the map $G$ is a Thom map, that is, a stratified submersion satisfying Thom's $a_G$ condition for every pair of strata (see \cite[Definitions~2.1 and~2.2]{GimenezConejero2023}), and its discriminant is a curve. This strengthens \cite[Theorem~5.7]{LeNunoSeade2025} under the additional assumption $\Sigma(f_0)=\{0\}$: the line $H_0$, which is part of the discriminant there, can be removed from the stratification of the target. This difference is essential to prove \cref{theo_main}, since the vector fields $\vec v_1,\vec v_2$ of \cref{section_3} do not preserve $H_0\times\Dbb_r$ (see \cref{fig:vfields}). We include the argument for the convenience of the reader, following also \cite[Theorem~2.6]{GimenezConejero2023}.

\begin{proposition}\label{prop_Thom}
Assume that $\Sigma(f_0)=\{0\}$ and that condition \eqref{eq:H} holds. Then, for small enough representatives:
\begin{enumerate}[(i)]
\item $\Sigma(G)$ is either $\{0\}$ or a curve, $G|_{\Sigma(G)}$ is finite, and $\Delta_G$ is either $\{0\}$ or a plane curve (see \cite[Theorem~5.7]{LeNunoSeade2025}).
\item There exist $\e_0,\nu_0>0$ and a Whitney stratification $\Scal_G$ of $X\cap\Bbb_{\e_0}\cap G^{-1}(\Dbb_{\nu_0})$ refining $\Scal$, with $\{0\}\in\Scal_G$, such that $G$ is a stratified submersion satisfying Thom's $a_G$ condition with respect to $\Scal_G$ and the stratification $\mathcal{T}=\big\{\Dbb_{\nu_0}\setminus\Delta_G,\ \Delta_G\setminus\{0\},\ \{0\}\big\}$, where $\Dbb_{\nu_0}\subset\C^2$ is the open disc of radius $\nu_0$.
\item For $0<\e<\e_0$ small enough and $U\subset\Dbb_{\nu_0}$ a small enough disc centered at $0$, the sphere $\mathbb{S}_\e$ is transverse to $G^{-1}(y)\cap A$ for every $y\in U$ and every stratum $A\in\Scal_G$. In particular, $G$ is transverse to small spheres off $\Delta_G$.
\end{enumerate}
\end{proposition}

\begin{proof}
(i) The set $\Sigma(G)$ is a closed analytic subset of $X$. If $Z$ is an irreducible component of $\Sigma(G)$ at $0$, then $Z\cap V_0=Z\cap f_0^{-1}(0)$ has dimension at least $\dim Z-1$ at $0$, and it is equal to $\{0\}$ by \eqref{eq:H2}; hence $\dim Z\leq1$. Moreover, $\Sigma(G)\cap G^{-1}(0)\subset\Sigma(G)\cap V_0=\{0\}$, so $G|_{\Sigma(G)}$ is finite and, by Remmert's finite mapping theorem, $\Delta_G=G\big(\Sigma(G)\big)$ is an analytic set of dimension at most one. After shrinking the representatives if necessary, $\Sigma(G)\setminus\{0\}$ and $\Delta_G\setminus\{0\}$ are smooth, each connected component of $\Sigma(G)\setminus\{0\}$ is contained in a stratum of $\Scal$, and $G$ restricts to a local diffeomorphism $\Sigma(G)\setminus\{0\}\to\Delta_G\setminus\{0\}$ (cf.\ \cite[Theorem~5.7]{LeNunoSeade2025}).

(ii) Let $\Scal_G$ consist of $\{0\}$, the connected components of $\Sigma(G)\setminus\{0\}$, and the nonempty sets $\big(S_\alpha\cap G^{-1}(T)\big)\setminus\Sigma(G)$, for $S_\alpha\in\Scal\setminus\{\{0\}\}$ and $T\in\mathcal{T}$. Since $G|_{S_\alpha\setminus\Sigma(G)}$ is a submersion, the latter sets are smooth, $G$ maps each of them submersively into the corresponding $T$, and they form a Whitney stratification of $X\setminus\Sigma(G)$ by \cite[Lemma~2.5]{GimenezConejero2023}. On the curves of $\Sigma(G)\setminus\{0\}$, $G$ is a local diffeomorphism onto $\Delta_G\setminus\{0\}$. Finally, shrinking the representative if necessary, we can assume that the Whitney conditions between the strata off $\Sigma(G)$ and the curves of $\Sigma(G)\setminus\{0\}$ are satisfied (see the proofs of \cite[Theorem~2.6]{GimenezConejero2023} and \cite[Theorem~5.7]{LeNunoSeade2025}).

Let us check Thom's condition. Let $A,B\in\Scal_G$ with $A\subset\overline B\setminus B$, and let $x_n\in B$ be a sequence converging to $x\in A$. 
If $A\subset\Sigma(G)$ the condition is trivial, because $\ker d_x(G|_A)=0$ since $A$ is either $\{0\}$ or a curve on which $G$ is a local diffeomorphism.
Otherwise, $A\subset S_\alpha\setminus\Sigma(G)$ for some $S_\alpha\in\Scal$, which implies that $B\subset S_\beta\setminus\Sigma(G)$ for some $S_\beta\in\Scal$. As $G|_{S_\beta}$ is a submersion along $B$, we have 
\[\ker d_{x_n}(G|_B)=\ker d_{x_n}(G|_{S_\beta})=T_{x_n}S_\beta\cap\ker d_{x_n}\widetilde G,\]
where $\widetilde G$ is an analytic extension of $G$ to $\C^N$. Passing to a subsequence, we may assume that the sequences $T_{x_n}S_\beta$ and $\ker d_{x_n}(G|_B)$ converge. By Whitney's condition (a), 
\[T_xS_\alpha\subset \lim T_{x_n}S_\beta.\] 
Since $d_x\widetilde G$ is surjective on $T_xS_\alpha$, it is also surjective on $\lim T_{x_n}S_\beta$, so the dimension of
\[\lim T_{x_n}S_\beta\cap\ker d_x\widetilde G\] is equal to the dimension of $\lim \ker d_{x_n}(G|_B)$ and they must be equal because the latter is contained in the former. Then, we obtain
\[
\lim \ker d_{x_n}(G|_B)=\lim T_{x_n}S_\beta \cap\ker d_x\widetilde G\supset T_xS_\alpha\cap\ker d_x\widetilde G\supset\ker d_x(G|_A),
\]
which is Thom's condition.

(iii) Let $\e>0$ be small enough so that $\mathbb{S}_{\e}$ is transverse to every stratum of $\Scal_G$. Away from the strata over $\Delta_G$, the statement is a standard consequence of the Thom condition. We include also the proof for any strata.

Assume that the claim fails for every disc $U$. Then there exist a stratum $B\in\Scal_G$ and points $x_n\in B\cap\mathbb{S}_\e$ with $G(x_n)\to 0$ such that $\ker d_{x_n}(G|_B)\subset T_{x_n}\mathbb{S}_\e$. Passing to a subsequence, $x_n$ converges to $x\in\mathbb{S}_\e\cap G^{-1}(0)$ and the limit of $\ker d_{x_n}(G|_B)$ exists. Since $x\neq0$ and $\Sigma(G)\cap G^{-1}(0)=\{0\}$, the point $x$ lies in a stratum $A\subset G^{-1}(0)$ not contained in $\Sigma(G)$, and $\ker d_x(G|_A)=T_xA$. By Thom's condition when $A\neq B$, or by smoothness when $A=B$, 
\[T_xA\subset \lim\ker d_{x_n}(G|_B) \subset T_x\mathbb{S}_\e,\]
 which contradicts the transversality of $A$ and $\mathbb{S}_\e$. 

 The last statement follows. Note that $\Delta_G\cap U$ does not depend on the representative of $G$ for $U$ small since $G|_{\Sigma(G)}$ is finite by (i).
\end{proof}

\begin{proof}[Proof of \cref{theo_main} for linear families]
Let $\e>0$ and $U$ be as in \cref{prop_Thom}\,(iii), with $U$ contained in the disc of \cref{remark_delta}, and let $r$ and $\vec v_1,\vec v_2$ be as in Steps~1 and~2 of \cref{section_3}. Set $M:=X\cap\overline{\Bbb}_\e\cap G^{-1}(U)$, stratified by the sets $A\cap\mathring{\Bbb}_\e$ and $A\cap\mathbb{S}_\e$, $A\in\Scal_G$. This is a Whitney stratification by \cite[Lemma~2.5]{GimenezConejero2023} applied to the function $\norm{x}^2$, since $\mathbb{S}_\e$ is transverse to the strata of $\Scal_G$ and does not meet the stratum $\{0\}$. We claim that
\[
\Gcal:M\times\Dbb_r\longrightarrow U\times\Dbb_r,\qquad(x,s)\longmapsto\big(G(x),s\big),
\]
is a proper Thom map with respect to this stratification (times $\Dbb_r$) and the stratification $\{T\times\Dbb_r\}_{T\in\mathcal T}$ of $U\times\Dbb_r$. 

Properness follows as in Step~3 of \cref{section_3}. 

By \cref{prop_Thom}\,(iii), the strata $A\cap\mathbb{S}_\e$ that are nonempty satisfy $\ker d_x(G|_A)\not\subset T_x\mathbb{S}_\e$, so $G$ is a stratified submersion.

\begingroup
Thom's condition for pairs of strata contained in
$\mathring{\Bbb}_\e$ follows from \cref{prop_Thom}\,(ii).
It remains to consider a limiting stratum contained in
$\mathbb{S}_\e$.

Let $x_n\in B'$ converge to $x\in A\cap\mathbb{S}_\e$, where
$A,B\in\Scal_G$ and
\[
B'=B\cap\mathring{\Bbb}_\e
\qquad\text{or}\qquad
B'=B\cap\mathbb{S}_\e.
\]
Assume that
\[
K'_n:=\ker d_{x_n}(G|_{B'})
\]
converges to a subspace $K'$.
Passing to a subsequence, we may also assume that
\[
K_n:=\ker d_{x_n}(G|_B)
\]
converges to a subspace $K$.
Set
\[
L:=\ker d_x(G|_A).
\]
By \cref{prop_Thom}\,(ii), or by smoothness when $A=B$, we have
$L\subset K$.
Moreover, \cref{prop_Thom}\,(iii) gives $L\not\subset T_x\mathbb{S}_\e$.
The kernel on the limiting boundary stratum is
\[
\ker d_x(G|_{A\cap\mathbb{S}_\e})=L\cap T_x\mathbb{S}_\e.
\]

If $B'=B\cap\mathring{\Bbb}_\e$, then $K'_n=K_n$, since $B'$ is
open in $B$. Hence
\[
K'=K\supset L\supset L\cap T_x\mathbb{S}_\e
=\ker d_x(G|_{A\cap\mathbb{S}_\e}),
\]
which proves Thom's condition in this case.

Suppose now that $B'=B\cap\mathbb{S}_\e$.
Then
\[
K'_n=K_n\cap T_{x_n}\mathbb{S}_\e,
\]
and therefore
\[
K'\subset K\cap T_x\mathbb{S}_\e.
\]
All dimensions in the following argument are real.
By \cref{prop_Thom}\,(iii), each $K_n$ is transverse to
$T_{x_n}\mathbb{S}_\e$, so
\[
\dim K'=\dim K-1.
\]
Since $L\subset K$ and $L\not\subset T_x\mathbb{S}_\e$, we also have
$K\not\subset T_x\mathbb{S}_\e$, and thus
\[
\dim(K\cap T_x\mathbb{S}_\e)=\dim K-1.
\]
The inclusion $K'\subset K\cap T_x\mathbb{S}_\e$ is consequently an equality.
It follows that
\[
K'=K\cap T_x\mathbb{S}_\e\supset L\cap T_x\mathbb{S}_\e
=\ker d_x(G|_{A\cap\mathbb{S}_\e}),
\]
which proves Thom's condition also in the boundary case.

Thus the restriction of $G$ to $M$, with the stated
stratification, satisfies Thom's condition.
Taking the product with $\Dbb_r$ preserves this condition.
Together with properness and stratified submersivity, already
established above, this proves that $\Gcal$ is a proper Thom map.

\par
\endgroup

\begingroup
By construction, the vector fields $\vec v_1,\vec v_2$ are
weakly controlled with respect to a suitable product tube
system for the target stratification: their parameter
components are constant, and they coincide with
$\partial/\partial s_j$ on a neighborhood of
$(\Delta_G\setminus\{0\})\times\Dbb_r$.
Since $\Gcal$ is a Thom map,
\cite[Theorems~3.5, 3.7 and~3.9]{GimenezConejero2023}
provide controlled stratified lifts $\vec u_1,\vec u_2$
which are locally integrable.
Properness allows their flows to be continued throughout
the parameter segments used in \cref{section_3}, as in
Step~3 of that section.

Applying the construction of Step~4 of \cref{section_3}
to these lifts gives a homeomorphism
\[
M\times\Dbb_r\longrightarrow M\times\Dbb_r,
\qquad (x,s)\longmapsto\big(h_s(x),s\big),
\]
with $h_0=\mathrm{id}$ and
\[
G\circ h_s=\psi_s^{-1}\circ G.
\]
By \cref{prop_vec_v}\,(ii), we obtain
\[
h_s(V_0\cap M)=V_s\cap M.
\]
Moreover, $h_s(0)=0$, since $\{0\}\times\Dbb_r$ is a
stratum preserved by the lifted flows.
Thus the maps $h_s$ induce the required homeomorphisms
of germs
\[
(X,0)\longrightarrow(X,0),
\qquad
(V_0,0)\longmapsto(V_s,0),
\]
proving \cref{theo_main} for linear families.

\par
\endgroup
\end{proof}

\begin{remark}\label{rem:othercondition}
Hypothesis \eqref{eq:H} is only used to show that $G$ is a Thom map (\cref{prop_Thom}), and the same proof works whenever $(f_s)$ is $\Delta$-regular, $\Delta_G$ is a curve, and $G$ is a Thom map for the stratification $\{U\setminus\Delta_G,\Delta_G\setminus\{0\},\{0\}\}$ of the target. See \cref{ex_2} below.
\end{remark}

We finish with two examples that illustrate the role of \eqref{eq:H} and the Thom condition.

\begin{example}\label{ex_1}
Let $X=\C^2$ with the stratification $\{\C^2\setminus\{0\},\{0\}\}$, $f_0=x^2+y^2$ and $\varphi=x^2$, so that $f_s=(1+s)x^2+y^2$ has a uniform critical set $\Sigma=\{0\}$. The critical locus of $\varphi$ is the line $\{x=0\}$, so $\varphi$ does not have an isolated critical point, and $\Sigma(G)=\{xy=0\}$. However, $\Sigma(G)\cap V_0=\{0\}$, so \eqref{eq:H} holds and \cref{theo_main} applies. In this case $\Delta_G=\{y_1=0\}\cup\{y_0=y_1\}$.
\end{example}

\begin{example}\label{ex_2}
Let $X=\C^3$ with the stratification $\{\C^3\setminus\{0\},\{0\}\}$, $f_0=x^2+y^2+z^2$ and $\varphi=x^2+y^2$. The family
\[
f_s=(1+s)(x^2+y^2)+z^2
\]
has a uniform critical set $\Sigma=\{0\}$, and it is clearly embeddedly trivial. However, $\Sigma(G)=\{z=0\}\cup\{x=y=0\}$ and
\[
\Sigma(G)\cap V_0=\{z=0,\ x^2+y^2=0\}
\]
is a pair of lines, so \eqref{eq:H} does not hold. Hence \eqref{eq:H} is not a necessary condition for embedded topological triviality, but it is easy to see that in this case \cref{rem:othercondition} applies and the same result holds.
\end{example}

\begin{proof}[Proof of \cref{cor_ICIS} for linear families]
Take $\Scal=\{X\setminus\{0\},\{0\}\}$. Since $\mu(f_s)$ is constant, the family has a uniform critical set $\Sigma=\{0\}$ (see \cite{Carvalho2019}), hence it is $\Delta$-regular by \cref{prop_MJA}. By \cref{rem_H}, the assumption on $X\cap f_0^{-1}(0)\cap\varphi^{-1}(0)$ is equivalent to \eqref{eq:H}, so the result follows from \cref{theo_main}.
\end{proof}

\section{\texorpdfstring{Separable families}{Separable families}}
\label{section_5}

We now deduce \cref{theo_A}, \cref{theo_main} and \cref{cor_ICIS} for a separable family $f_s=f_0+\rho(s)\varphi$, with $\rho\not\equiv0$, from the linear case proved in \cref{section_3,section_4}. Let
\[
\ell_t:=f_0+t\,\varphi,\qquad t\in\C,
\]
be the linear family with the same $f_0$ and $\varphi$. Then $f_s=\ell_{\rho(s)}$, that is, $(f_s)$ is obtained from $(\ell_t)$ by the change of parameter $t=\rho(s)$. Both families have the same map $G$, so $\Delta$-regularity, the transversality of \cref{def:ML}, condition \eqref{eq:H} and the set $N$ are the same for both. The remaining hypotheses are conditions on the members of the family, which are handled as follows.

\begin{lemma}\label{lem_separable}
Let $P$ be a property of function-germs on $(X,0)$. Then $f_s$ has $P$ for every $s$ in a neighborhood of $0$ if and only if $\ell_t$ has $P$ for every $t$ in a neighborhood of $0$.
\end{lemma}

\begin{proof}
Since $\rho(0)=0$ and $\rho\not\equiv0$, the map $\rho$ is continuous and, by the open mapping theorem, it sends neighborhoods of $0$ to neighborhoods of $0$. Hence, if $\ell_t$ has $P$ for $t\in D'$, then $f_s=\ell_{\rho(s)}$ has $P$ for $s\in\rho^{-1}(D')$, and if $f_s$ has $P$ for $s\in D$, then $\ell_t$ has $P$ for $t\in\rho(D)$.
\end{proof}

Applying \cref{lem_separable} to the properties $\Sigma(\,\cdot\,)\cap\Bbb_\e=\Sigma(f_0)\cap\Bbb_\e$, $\mu(\,\cdot\,)=\mu(f_0)$ and having an isolated singularity, we see that $(f_s)$ has a uniform critical set, has constant Milnor number, or consists of isolated singularity function-germs if and only if $(\ell_t)$ does. In particular, \cref{prop_MJA} for $(f_s)$ follows from the linear case \cite{Menegon2023}.

\begin{proof}[Proof of \cref{theo_A}, \cref{theo_main} and \cref{cor_ICIS} for separable families]
By the above, $(\ell_t)$ satisfies the hypotheses of the statement. Hence, by the linear case, there are $M$, a disc $D'$ centered at $0$ and a homeomorphism $(x,t)\mapsto\big(h_t(x),t\big)$ of $M'\times D'$ satisfying the statement for $(\ell_t)$, where $M'=M\setminus N$ for \cref{theo_A} and $M'=M$ otherwise; in particular, $h_0=\mathrm{id}$ and $h_t(V_0\cap M')=\ell_t^{-1}(0)\cap M'$. Let $D$ be a disc centered at $0$ with $\rho(D)\subset D'$. Then
\[
M'\times D\to M'\times D,\qquad(x,s)\mapsto\big(h_{\rho(s)}(x),s\big),
\]
is a homeomorphism, since it is continuous and so is its inverse $(x,s)\mapsto\big(h_{\rho(s)}^{-1}(x),s\big)$. Moreover, $h_{\rho(0)}=h_0=\mathrm{id}$,
\[
h_{\rho(s)}(V_0\cap M')=\ell_{\rho(s)}^{-1}(0)\cap M'=V_s\cap M',
\]
and each $h_{\rho(s)}$ preserves the strata (and fixes $0$, for \cref{theo_main}) because each $h_t$ does. Hence the maps $h_{\rho(s)}$, $s\in D$, satisfy the statement for $(f_s)$.
\end{proof}

\subsection*{AI assistance}

A large language model was used to assist in revising
mathematical arguments and their presentation. The authors
take full responsibility for the final manuscript.



\end{document}